\documentclass[conference]{IEEEtran}
\usepackage{amssymb}
\usepackage{amsmath}
\usepackage[bookmarks=false]{hyperref}
\usepackage{graphicx}            
\usepackage{amsfonts}              
\usepackage{amsthm}             
\usepackage{algorithm}
\usepackage{algpseudocode}
\usepackage{cite}

\hypersetup{
colorlinks = false,
citecolor = blue,
urlcolor = blue,
linkcolor = blue
}

\newtheorem{thm}{Theorem}

\begin{document}

\title{ Spectrum-Revealing Cholesky Factorization for Kernel Methods }

\author{\IEEEauthorblockN{Jianwei Xiao}
\IEEEauthorblockA{Department of Mathematics\\
UC Berkeley, Berkeley, CA\\
\href{mailto:jwxiao@berkeley.edu}{jwxiao@berkeley.edu}}
\and
\IEEEauthorblockN{Ming Gu}
\IEEEauthorblockA{Department of Mathematics\\
UC Berkeley, Berkeley, CA\\\
\href{mailto:mgu@math.berkeley.edu}{mgu@math.berkeley.edu}}
}

\maketitle

\begin{abstract}
Kernel methods represent some of the most popular machine learning tools for data analysis. Since exact kernel methods can be prohibitively expensive for large problems, reliable low-rank matrix approximations and high-performance implementations have become indispensable for practical applications of kernel methods. In this work, we introduce spectrum-revealing Cholesky factorization, a reliable low-rank matrix factorization, for kernel matrix approximation. 
We also develop an efficient and effective randomized algorithm for computing this factorization. 
Our numerical experiments demonstrate that this algorithm is as effective as other Cholesky factorization based kernel methods on machine learning problems, but significantly faster.
\end{abstract}

\IEEEpeerreviewmaketitle 

\section{Introduction}\label{section_Introduction}

\subsection{Kernel methods and their practical performance}\label{subsection_Kernel methods and their practical performance}

Kernel methods have become an increasingly popular tool for machine learning tasks such as classification, prediction and clustering, with diverse applications including information extraction and handwriting recognition. 

Kernel methods owe their name to the use of kernel functions, which enable them to operate in a high-dimensional feature space. One critical drawback of kernel methods is their inability to solve very large-scale learning problems owing to their high cost. Given $n$ data points, the kernel matrix $K$ is of size $n \times n$, which implies a computational complexity of at least $O(n^2)$. More importantly, most kernel methods require matrix inversions or eigen-decompositions as their computational core, leading to complexities as high as $O(n^3)$. 

There are two major approaches to significantly improving the performance of kernel methods. Firstly, the aforementioned complexities can be reduced by approximating the kernel matrices with their low-rank approximations, and this is one of the major reasons for the practical success of kernel methods. The kernel matrix $K$ can be approximated in the form 
\begin{equation}\label{Eqn:KLLT}
K \approx LL^T, \quad \mbox{where} \quad L\in \mathbb{R}^{n \times k}, 
\end{equation}
and approximate rank $k$ is chosen so that $k \ll n$. Moreover, it is often possible to reformulate kernel methods to make use of $L$ instead of $K$. This result in learning methods of much lowered computational complexity of $O(k^3+k^2n)$ \cite{fine2002efficient, bach2003kernel}.

Another approach is to develop highly-tuned software libraries for machine learning algorithms for kernel methods \cite{chang2011libsvm,pedregosa2011scikit}. By reorganizing their internal computations to take advantage of high performance linear algebra packages such as LAPACK \cite{anderson1999lapack}, significant practical speedups can be produced without major mathematical changes to the algorithms. 

\subsection{Cholesky factorization based low-rank approximation}\label{subsection_Cholesky factorization based low-rank approximation} 

One of the most popular methods to obtain a low-rank approximation of a kernel matrix in the form (\ref{Eqn:KLLT}) is based on the Cholesky factorization. Such approximations have been used in many areas: SVM training \cite{fine2002efficient}, kernel independent component analysis \cite{bach2003kernel} and predictive low-rank decomposition for kernel methods \cite{bach2005predictive}. The essential part in finding a low-rank approximation of kernel matrix using Cholesky factorization is to find representative training samples, which is equivalent to doing Cholesky factorization of kernel matrix with certain pivoting strategy. Diagonal pivoting is commonly used in these Cholesky factorization based algorithms.

However, there are two major well-known drawbacks of diagonal pivoting. Firstly, pivots are computed one at a time, which results in mostly level-2 BLAS operations, much less efficient than level-3 BLAS operations. Secondly, the pivots chosen by diagonal pivoting may occasionally fail to produce a reliable low-rank approximation to the kernel matrix. 

\subsection{Randomized Cholesky for reliable low-rank approximation}\label{subsection_Randomized Cholesky for reliable low-rank approximation}

In recent works, randomization has emerged as an especially potent tool for large-scale data analysis. Reliable and efficient randomized algorithms have been successfully developed for low-rank approximation of matrices \cite{gu2015subspace}, sketching problems \cite{demmel2015communication}, and fast solution to the least squares problem \cite{drineas2011faster}.

In this work, we develop a randomized blocked Cholesky factorization algorithm for reliable low-rank approximation of the kernel matrix. This algorithm is run in two stages. In the first stage, we first use randomization to project the kernel matrix into a smaller one; we then find pivots on the smaller matrix, and finally we apply these pivots to the kernel matrix in a block form to fully take advantage of level-3 BLAS. We repeat this process until we reach approximate rank $k$. Despite the randomness, this approach works really well in practice. 

However, occasionally the approach above may not lead to reliable low-rank approximations. To guard against this possibility, in the second stage we further employ a novel follow-up pivoting strategy that simultaneously ensures a reliable low-rank approximation and separates the representative training samples from those that are nearly linearly dependent on them. This separation feature is of significant interest in its own right \cite{wold1984collinearity}. We will establish strong singular value and matrix error bounds to demonstrate the effectiveness of this pivoting strategy. 

Although the computational complexity of our algorithm is no longer linear in $n$, the implementation of our algorithm is still faster than other Cholesky factorization based algorithms. The main reason is that run time of an algorithm is not only dependent on arithmetic cost but also dependent on communication cost, which represents the required time to move data. The communication cost often dominates arithmetic cost. Level-3 BLAS have significantly lower communication cost than level-2 BLAS. As a blocked algorithm, our novel method fully utilizes level-3 BLAS. As our experiments demonstrate, while our method is at least as reliable as other Cholesky factorization based algorithms in all applications, it is much faster for large scale low-rank approximations. 

\subsection{Our Contributions}\label{subsection_Our Contributions}

\textbf{Spectrum-revealing Cholesky factorization:} We demonstrate the existence of the Spectrum-revealing Cholesky factorization (SRCH), and develop strong singular value and matrix error bounds for SRCH. Our analysis shows that SRCH provides a highly reliable low-rank approximation to the kernel matrix for any given approximate rank $k$. 

\textbf{A randomized algorithm for computing an SRCH:} Unlike existing Cholesky factorization algorithms, this randomized algorithm can efficiently and correctly compute an SRCH. It is especially suitable to obtain quality low-rank approximations for matrices with fast decaying singular-value spectra, which are ubiquitous in kernel matrices in machine learning.

\textbf{Empirical validation:} We compare our method with other Cholesky factorization based algorithms in two different applications: a prediction problem and the Gaussian process. All of these methods show similar effectiveness, but our method is significantly faster.

\subsection{Paper Summary}\label{subsection_Paper Summary}

In section \ref{section_The Setup and Background} we briefly introduce previous work of pivoted Cholesky. In section \ref{subsection_A randomized blocked left-looking Cholesky factorization} we develop a randomized blocked left-looking algorithm to compute a pivoted Cholesky factorization without explicitly updating the Schur complement. In section \ref{subsection_Spectrum-revealing Cholesky factorization} we define and discuss properties of a Spectrum-revealing Cholesky factorization (SRCH); and develop an efficient modification to the algorithm in section \ref{subsection_A randomized blocked left-looking Cholesky factorization} to reliably compute an SRCH. In section \ref{section_Numerical experiments} we compare our algorithm and other alternatives in different applications. In section \ref{section_Conclusion} we conclude this work. 

\section{The Setup and Background}\label{section_The Setup and Background}

\subsection{Notation}\label{subsection_Notation}

$\sigma_j(A)$ denotes the $j$th largest singular value of $A$. If $A$ is real symmetric, $\lambda_j(A)$ denotes its $j$th largest eigenvalue. In this work we follow MATLAB notation.

\subsection{Diagonal pivoted Cholesky factorization}\label{subsection_Diagonal pivoted Cholesky factorization}

Diagonal pivoting is the most popular pivoting strategy in computing a Cholesky factorization for low-rank approximation. This strategy chooses the largest diagonal entry as the pivot at each pivoting step. \cite{lucas2004lapack} developed {\tt {\tt DPSTRF.f}} for the diagonal pivoted Cholesky factorization.

There are two major problems with {\tt {\tt DPSTRF.f}}. Firstly, most of its work is in updating the Schur complement. However, the Schur complement on exit is typically discarded in a low-rank matrix approximation, meaning most of this work is unnecessary if $k \ll n$. Secondly, diagonal pivoting is a greedy strategy for computing a low-rank approximation by pivoting to the largest diagonal entry. There are well-known classes of matrices for which this strategy fails to compute a reliable low-rank approximation \cite{higham2002accuracy}. \cite{gu2004strong} provides an algorithm that can always compute a reliable low-rank approximation by doing suitable swaps after obtaining a partial Cholesky factorization with diagonal pivoting, but their algorithm is not very efficient.

\section{New algorithms and main results}\label{section_New algorithms and main results}

\subsection{A randomized blocked left-looking Cholesky factorization}\label{subsection_A randomized blocked left-looking Cholesky factorization}

The Cholesky factorization can be computed in a number of different, but mathematically equivalent, variants. \hyperref[Left Looking Blocked Cholesky factorization]{Algorithm \ref{Left Looking Blocked Cholesky factorization}} 
is a left-looking variant that computes the full Cholesky factorization without directly updating the Schur complement. For a symmetric positive definite $X$, ${\bf chol}(X)$ is the Cholesky factor such that $\left({\bf chol}(X)\right) \left({\bf chol}(X)\right)^T = X$.

\begin{algorithm}
\caption{Left Looking Blocked Cholesky factorization}\label{Left Looking Blocked Cholesky factorization}
\begin{algorithmic}
\State $\textbf{Inputs:}$  $\;\;\;\;$ Positive semidefinite $A\in{}\mathbb{R}^{n\times{}n}$; block size $b$. 
\State $\textbf{Outputs:\;\,\,}$   $L\stackrel{def}{=}$ lower triangular part of $A$.  
\For{$\mbox{j=1:b:n}$}
\begin{eqnarray}
\overline{b}\!\!\!\!&\;\;=&\!\!\!\!{\bf min}(\mbox{b,n-j+1})  \nonumber \\
A(\mbox{j:n,j:j+$\overline{b}$-1})\!\!\!\!&\mbox{--$\;$=}&\!\!\!\!A(\mbox{j:n,1:j-1)A(j:j+$\overline{b}$-1,1:j-1})^T \label{Eqn:LUpdate2} \\
A(\mbox{j:j+$\overline{b}$-1,j:j+$\overline{b}$-1})\!\!\!\!&\,\;=&\!\!\!\!{\bf chol}(A(\mbox{j:j+$\overline{b}$-1,j:j+$\overline{b}$-1})) \nonumber \\
A(\mbox{j+$\overline{b}$:n,j:j+$\overline{b}$-1})\!\!\!\!&\,\;=&\!\!\!\!A(\mbox{j+$\overline{b}$:n,j:j+$\overline{b}$-1})  A(\mbox{j:j+$\overline{b}$-1,j:j+$\overline{b}$-1})^{-T} \nonumber 
\end{eqnarray}
\EndFor
\end{algorithmic}
\end{algorithm}

Most of the work in \hyperref[Left Looking Blocked Cholesky factorization]{Algorithm \ref{Left Looking Blocked Cholesky factorization}} is in updating the matrix $A(\mbox{j:n,j:j+$\overline{b}$-1})$ in equation (\ref{Eqn:LUpdate2}). This work starts small and increases linearly with $j$. Thus, \hyperref[Left Looking Blocked Cholesky factorization]{Algorithm \ref{Left Looking Blocked Cholesky factorization}} would be much faster than {\tt {\tt DPSTRF.f}} if we restricted $j \leq k$ for some $k \ll n$. But such restriction, without the benefit of pivoting, may not lead to a very meaningful approximation of $A$. 

Based on this consideration, we now introduce a novel pivoting strategy into \hyperref[Left Looking Blocked Cholesky factorization]{Algorithm \ref{Left Looking Blocked Cholesky factorization}}. For given $p \geq b$ and $A \in \mathbb{R}^{n \times n}$, we draw a random matrix  $\Omega \in \mathcal{N}(0,1)^{p \times n}$ and compute a random projection $B=\Omega A$, which is significantly smaller than $A$ in row dimension if $p \ll n$. We compute a partial QR factorization with column pivoting (QRCP) on $B$ to obtain $b$ column pivots and apply them as $b$ diagonal pivots on $A$. Intuitively, good pivots for $B$ should also be good pivots for $A$. For this strategy to work, we need to compute a random projection for the Schur complement $A(\mbox{j+$\overline{b}$:n,j+$\overline{b}$:n})$ for each $j$, without explicitly computing $A(\mbox{j+$\overline{b}$:n,j+$\overline{b}$:n})$. Remarkably, such a random projection can indeed be quickly computed via an updating formula. \hyperref[algorithm_SRCH]{Algorithm \ref{algorithm_SRCH}} computes a partial Cholesky factorization, with diagonal pivots chosen by partial QRCP on successive random projections.

\begin{algorithm}
\caption{Randomized blocked left-looking Cholesky factorization}\label{algorithm_SRCH}
\begin{algorithmic}
\State $\textbf{Inputs:}$  $\;\;\;\;$ Positive semidefinite $A\in{}\mathbb{R}^{n\times{}n}$; 
\State $\quad \quad \quad \quad \;\; $ block size $b$; over-sampling size $p \geq b$; 
\State $\quad \quad \quad \quad \;\; $ approximate rank $k \ll n$. 
\State $\textbf{Outputs:\;\,\,}$ Permutation vector ${\Pi}\in{}\mathbb{R}^n$; 
\State $\quad \quad \quad \quad \;\; $  $L\stackrel{def}{=}$ lower triangular part of $A(\mbox{1:n,1:k})$.  
\State $\textbf{Initialize:\,\,}$ {\bf generate} random $\Omega \in \mathcal{N}(0,1)^{p\times{}n}$; 
\State $\quad \quad \quad \quad \;\; $  {\bf compute} $B=\Omega{}A$, $\Pi = (\mbox{1:n}) \in \mathbb{R}^n$.
\For{$\mbox{j=1:b:k}$}
\State $\overline{b}={\bf min}\mbox{(b,k-j+1)}$
\State {\bf compute} partial QRCP on $B(\mbox{:,j:n})$ to obtain $\overline{b}$ 
\State $\quad  $ column pivots $(\mbox{j}_1,\mbox{j}_2,\dots,\mbox{j}_{\overline{b}})$
\State {\bf swap} $(\mbox{$\mbox{j}_1$,$\mbox{j}_2$,$\dots$,$\mbox{j}_{\overline{b}}$})$ and $(\mbox{$\mbox{j}$,$\mbox{j+1}$,\dots,$\mbox{j+}\overline{b}\mbox{-1}$})$ columns 
\State $\quad $  in $B, \Omega$, and corresponding entries in $\Pi$
\State {\bf swap} corresponding rows and columns in $A$
\begin{eqnarray}
A(\mbox{j:n,j:j+$\overline{b}$-1}) \!\!\!\!&\mbox{--$\;$=}&\!\!\!\! A(\mbox{j:n,1:j-1)A(j:j+$\overline{b}$-1,1:j-1})^T  \label{Eqn:LUpdate3}\\
A(\mbox{j:j+$\overline{b}$-1,j:j+$\overline{b}$-1}) \!\!\!\!&\,=&\!\!\!\! {\bf chol}(A(\mbox{j:j+$\overline{b}$-1,j:j+$\overline{b}$-1})) \nonumber \\
A(\mbox{j+$\overline{b}$:n,j:j+$\overline{b}$-1}) \!\!\!\!&\,=&\!\!\!\! A(\mbox{j+$\overline{b}$:n,j:j+$\overline{b}$-1})  A(\mbox{j:j+$\overline{b}$-1,j:j+$\overline{b}$-1})^{-T} \nonumber 
\end{eqnarray}
\If{$j+\overline{b}-1 < k$}
\State $\hspace*{-0.1in} B(\mbox{:,j+$\overline{b}$:n}) \,\;\mbox{--$\;$=}\;\, \Omega(\mbox{:,j:n}) A(\mbox{j:n,j:j+$\overline{b}$-1}) A(\mbox{j+$\overline{b}$:n,j:j+$\overline{b}$-1})^T$
\EndIf
\EndFor
\end{algorithmic}
\end{algorithm}

\begin{figure}
\begin{center}
\includegraphics[height=4cm,width=8cm]{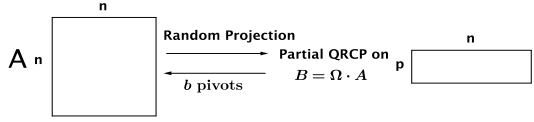}
\end{center}
\caption{Use random projection to find $b$ pivots in each iteration.}\label{figure_Use random projection to find a block of pivots in each iteration}
\end{figure}

\hyperref[figure_Use random projection to find a block of pivots in each iteration]{Fig. \ref{figure_Use random projection to find a block of pivots in each iteration}} illustrates the main ideas in \hyperref[algorithm_SRCH]{Algorithm \ref{algorithm_SRCH}}. We recursively find $b$ pivots on $B$ and apply these pivots on $A$.

Since \hyperref[algorithm_SRCH]{Algorithm \ref{algorithm_SRCH}} is a left-looking algorithm, it's much more efficient than {\tt {\tt DPSTRF.f}} when $k \ll n$. Numerical experiments also suggest that \hyperref[algorithm_SRCH]{Algorithm \ref{algorithm_SRCH}} typically computes a better low-rank approximation than {\tt {\tt DPSTRF.f}}.

\noindent {\bf Updating formula for $B$:} The formula for successively computing $B(\mbox{:,j+$\overline{b}$:n})$ for increasing $\mbox{j}$ in \hyperref[algorithm_SRCH]{Algorithm \ref{algorithm_SRCH}} is what makes \hyperref[algorithm_SRCH]{Algorithm \ref{algorithm_SRCH}} so efficient. To derive it, we first compute $B = \Omega \, A$. \hyperref[algorithm_SRCH]{Algorithm \ref{algorithm_SRCH}} then computes $b$ pivots based on a partial QRCP on $B$ and performs a block Cholesky step. To continue, \hyperref[algorithm_SRCH]{Algorithm \ref{algorithm_SRCH}} need to compute a random projection on the corresponding Schur complement, which it does not compute directly. We can re-use the initial random matrix $\Omega$ and the corresponding random projection $B$ to compute a special random projection for the Schur complement.

After the necessary row and column swaps and the block Cholesky step, we can write the matrices $A$, and $\Omega$ as ${\displaystyle \left(
\begin{array}{cc}
L_{11} &  \\
L_{21} & I \\
\end{array}
\right)\left(
\begin{array}{cc}
L_{11}^T & L_{21}^T \\
& S_2 \\
\end{array}
\right)}$,  and ${\displaystyle \left(
\begin{array}{cc}
\Omega_1 & \Omega_2 \\
\end{array}
\right)}
$, respectively. The column swapped $B$ can be written as
\[ {\displaystyle \left(
\begin{array}{cc}
B_1 & B_2 \\
\end{array}
\right) = \left(
\begin{array}{cc}
\Omega_1 & \Omega_2 \\
\end{array}
\right) \; \left(\begin{array}{cc}
L_{11} &  \\
L_{21} & I \\
\end{array}
\right)\left(
\begin{array}{cc}
L_{11}^T & L_{21}^T \\
& S_2 \\
\end{array}
\right),}\]
which, in turn, implies a special random projection formula
\begin{equation}\label{equation_updating_formula_B}
\Omega_{2}\;S_2=B_{2}-\left(
\begin{array}{cc}
\Omega_1 & \Omega_2 \\
\end{array}
\right)\left(\begin{array}{c}
L_{11}  \\
L_{21}  \\
\end{array}
\right)L_{21}^T.
\end{equation}
We can compute the random projection $\Omega_2 S_2$ for $S_2$ via the right hand side expression in equation (\ref{equation_updating_formula_B}), with $\Omega_2$ as the random matrix. Generalizing this consideration for all $\mbox{j}$ results in the formula for computing $B(\mbox{:,j+$\overline{b}$:n})$ in \hyperref[algorithm_SRCH]{Algorithm \ref{algorithm_SRCH}}. 

\noindent {\bf Complexity analysis:} The most work of \hyperref[algorithm_SRCH]{Algorithm \ref{algorithm_SRCH}} are computing the initial random projection matrix $B = \Omega A$ and updating pivoted block $A(\mbox{j:n,j+$\overline{b}$-1})$ for each $j$. Computing $B$ requires $O(n^2p)$ operations and updating all pivoted blocks requires $O(nk^2)$ operations, leading to the overall complexity of $O(n^2p+nk^2)$ operations. Note that if both $p\ll n$ and $k\ll n$, then the complexity is $O(n^2p)$, i.e., the dominant computation is in the overhead -- computing $B = \Omega A$.

\subsection{Spectrum-revealing Cholesky factorization}\label{subsection_Spectrum-revealing Cholesky factorization}
Greedy pivoting strategies in {\tt {\tt DPSTRF.f}} and \hyperref[algorithm_SRCH]{Algorithm \ref{algorithm_SRCH}} typically compute good quality low-rank approximations, but not always. Below we first discuss what low-rank approximations are possible based on diagonal pivoting alone, and then develop a swap strategy to modify the partial Cholesky factorization computed from \hyperref[algorithm_SRCH]{Algorithm \ref{algorithm_SRCH}} to ensure such an approximation. With a slight abuse of notation, $\Pi$ in \hyperref[theorem_swap_strategy]{Theorem \ref{theorem_swap_strategy}} denotes a permutation matrix. Recall that for any matrix $X$, $\left|\left|X\right|\right|_{2,1}$ is equal to the largest of the column norms of $X$, and  $\left|\left|X\right|\right|_{\max}$ is equal to the largest entry of $X$ in absolute value. 

\begin{thm}\label{theorem_swap_strategy} Let $A \in{}\mathbb{R}^{n\times{}n}$ be symmetric positive definite with a partial Cholesky factorization for $ k < n$: 
\begin{equation}
\Pi^TA\Pi=LL^T+\left(
\begin{array}{cc}
0 & 0 \\
0 & S \\
\end{array}
\right) \label{Eqn:SRCH}
\end{equation}
\[ \mbox{with } \quad L = \left(
\begin{array}{c}
L_{11} \\
\ell^T \\
L_{21} \\
\end{array}
\right), \quad \mbox{and} \quad S = \left(
\begin{array}{cc}
\alpha & s^T \\
s & \widehat{S} \\
\end{array}
\right), \]
where $L_{11}\in{}\mathbb{R}^{k\times{}k}, \ell\in{}\mathbb{R}^{k\times{}1}, L_{21}\in{}\mathbb{R}^{(n-k-1)\times{}k}, \alpha{}\in{}\mathbb{R}, s\in{}\mathbb{R}^{(n-k-1)\times{}1}$, and $ \widehat{S}\in{}\mathbb{R}^{(n-k-1)\times{}(n-k-1)}$.

Assume that $\alpha =  \left|\left|S\right|\right|_{\max}$ and that for a given $ g >1$,

\begin{align}\label{equation_rank_revealing_condition}
\frac{1}{\sqrt{\alpha}}&\geq{}\frac{1}{\sqrt{g}}\left|\left|\left(
\begin{array}{cc}
L_{11} &  \\
\ell^T & \sqrt{\alpha} \\
\end{array}
\right)^{-1}
\right|\right|_{2,1}, 
\end{align}
then there exists a $\tau \leq g(n-k)(k+1)$ such that for $1 \leq j \leq k$ 
\begin{eqnarray}
{\displaystyle \lambda_{k+1}(A) \leq \left\|\Pi^TA\Pi-L\,L^T\right\|_2 \leq \tau\; \lambda_{k+1}(A),}  \label{Eqn:thm2matnorm}  && \\
\displaystyle \lambda_{j}(A)  \geq   \sigma_{j}^{2}(L)}\geq{\displaystyle\frac{\lambda_{j}(A)}{1+\tau\; {\bf min}{}\left\{1,(1+\tau){\displaystyle \frac{\lambda_{k+1}(A)}{\lambda_{j}(A)}}\right\}}. \label{Eqn:thm2sv}  &&
\end{eqnarray}
\qed
\end{thm}

If conditions on $\alpha$ hold, then \hyperref[theorem_swap_strategy]{Theorem \ref{theorem_swap_strategy}} asserts that the matrix approximation error is at most a factor of $\tau$ away from being optimal in $2$-norm; and all the singular values of $L$ are at most a factor of $\sqrt{1 + \tau}$ away from being optimal. In addition, for the largest  singular values of $A$ where ${\displaystyle \frac{\lambda_{k+1}(A)}{\lambda_j(A)} \ll 1}$, the corresponding approximate singular values $\sigma_j(L)$ are very close to $\sqrt{\lambda_j(A)}$, the best possible for any rank $k$ approximation. 

A partial Cholesky factorization of the form (\ref{Eqn:SRCH}) is said to be {\em Spectrum-revealing} if it satisfies the conditions on $\alpha$ in  \hyperref[theorem_swap_strategy]{Theorem \ref{theorem_swap_strategy}}. The singular value lower bound in relation (\ref{Eqn:thm2sv}) represents a unique feature in a Spectrum-revealing Cholesky (SRCH) factorization.

{\tt {\tt DPSTRF.f}} ensures the condition $\alpha =  \left|\left|S\right|\right|_{\max}$ by performing diagonal pivoting, but may not  satisfy condition (\ref{equation_rank_revealing_condition}), leading to potentially poor approximations. 

\hyperref[algorithm_swap_strategy]{Algorithm \ref{algorithm_swap_strategy}} is a randomized algorithm that efficiently computes a Spectrum-revealing Cholesky factorization. It initializes the permutation with \hyperref[algorithm_SRCH]{Algorithm \ref{algorithm_SRCH}}. If \hyperref[algorithm_SRCH]{Algorithm \ref{algorithm_SRCH}} fails condition  (\ref{equation_rank_revealing_condition}),  \hyperref[algorithm_swap_strategy]{Algorithm \ref{algorithm_swap_strategy}} makes additional randomized column and row swaps to ensure it. In the algorithm, we denote $\widehat{L}=\left(
\begin{array}{cc}
L_{11} & \\
\ell^T &\sqrt{\alpha} \\
\end{array}
\right)
$.

\begin{algorithm}
\caption{A randomized algorithm to compute an SRCH}\label{algorithm_swap_strategy}
\begin{algorithmic}
\State $\textbf{Inputs:}$  $\;\;\;\;$ Positive semidefinite $A\in{}\mathbb{R}^{n\times{}n}$; 
\State $\quad \quad \quad \quad \;\; $ block size $b$; over-sampling size $p \geq b$; 
\State $\quad \quad \quad \quad \;\; $ parameter $g>1$; approximate rank $k \ll n$. 
\State $\textbf{Outputs:\;\,\,}$ Permutation vector ${\Pi}\in{}\mathbb{R}^{n}$; matrix $L$. 
\State $\textbf{Initialize:\,\,}$ {\bf compute} ${\Pi}$ and $L$ with \hyperref[algorithm_SRCH]{Algorithm \ref{algorithm_SRCH}} 
\State {\bf Generate} random $\Omega \in\mathcal{N}(0,1)^{d\times{}(k+1)}$ 
\State {\bf compute} $\alpha={\bf max}{}(\text{diag(Schur complement)})$
\While{${\displaystyle \frac{1}{\sqrt{\alpha}}<\frac{1}{\sqrt{gd}}\left\|\Omega{}\widehat{L}^{-1}\right\|_{2,1}}$}
\State $\imath ={\bf argmax}_{1\leq{}i\leq{}k+1}\{i\text{th column norm of }\Omega{}\widehat{L}^{-1}\}$
\State {\bf swap} $\imath$-th and $(k+1)$-st columns and rows of $A$
\State {\bf swap} $\imath$-th and $(k+1)$-st entries in $\Pi$
\State {\bf Givens-rotate} $L$ back into lower-triangular form.  
\State {\bf compute} $\alpha={\bf max}{}(\text{diag(Schur complement)})$
\EndWhile
\end{algorithmic}
\end{algorithm}

Remarks: \\
(1) We look through the diagonal of the Schur complement $S$ to find $\alpha$, thereby avoid computing $S$ itself. \\
(2) The while loop in \hyperref[algorithm_swap_strategy]{Algorithm \ref{algorithm_swap_strategy}} will eventually stop, after a finite number of swaps (bounded above by $O(n)$), leading to a permutation that satisfies conditions (\ref{equation_rank_revealing_condition}). However, at most a few swaps are enough in practice. \\
(3) Each swap will make the $\imath$-th row out of the lower-triangular form. A round robin rotation is applied to the rows of $L$ and a quick sequence of Givens rotations are right multiplied to $L$ to restore its lower-triangular form. These Givens rotations are orthogonal and will cancel themselves out in the matrix product $L L^T$. \\
(4) In some practical applications where more accurate singular values are desirable, one can compute an SRCH for a rank $\widehat{k}$ approximation $\widehat{L} \widehat{L}^T$ with $\widehat{k} > k$ and then SVD-truncate the matrix $\widehat{L}$ into a rank $k$ matrix $L$. This will lead to a rank-$k$ approximation $L L^T$ that satisfies
\small
\begin{eqnarray*}
\lambda_{k+1}(A) \leq \left\|\Pi^T A \Pi - LL^T \right\|_2 \leq \lambda_{k+1}(A) + \widehat{\tau} \lambda_{\widehat{k}+1}(A), \\
\lambda_j(A) \geq \sigma_j^2(L) \geq \frac{\lambda_j(A)}{1+\widehat{\tau} {\bf min}\left\{1,(1+\widehat{\tau})\frac{\lambda_{\widehat{k}+1}(A)}{\lambda_j(A)}\right\}}, 
\end{eqnarray*}
\normalsize
for $1 \leq j \leq k$ and scalar $\widehat{\tau} \leq g(n-\widehat{k})(\widehat{k}+1)$.
Especially for rapidly decaying singular values, i.e., $\lambda_{\widehat{k}+1}(A) \ll \lambda_{k+1}(A)$, these bounds make $L L^T$ almost indistinguishable from the best possible, the SVD-truncated rank-$k$ approximation.

\noindent{\bf Complexity analysis:} In addition to initialization, \hyperref[algorithm_swap_strategy]{Algorithm \ref{algorithm_swap_strategy}} needs to repeatedly compute the diagonal of the Schur complement, which can be done in $O(k(n-k))$ operations. We need to swap the column with largest diagonal entry to the leading column and update the pivoted column, which needs a matrix-vector multiplication, costing another $O(k(n-k))$ operations. Then we need to compute $\Omega{}\widehat{L}^{-1}$ and the corresponding column norms, costing $O(dk^2)$ operations. The algorithm stops if the while loop condition fails. Otherwise we need to swap the column with largest column norm in $\Omega{}\widehat{L}^{-1}$ and the leading column in the Schur complement. The Givens rotations needed to restore lower-triangular form in $L$ cost $O(nk)$ operations. In total, the cost of performing one swap in \hyperref[algorithm_swap_strategy]{Algorithm \ref{algorithm_swap_strategy}} is $O(nk+dk^2)$ operations. Assuming $k\ll{}n$, this cost becomes $O(nk)$, which is negligible compared to $O(n^2p)$, the cost of initialization (\hyperref[algorithm_SRCH]{Algorithm \ref{algorithm_SRCH}}). 

\section{Numerical experiments}\label{section_Numerical experiments}

Data and source code are available at following URL \footnote{\href{https://math.berkeley.edu/~jwxiao/}{\url{https://math.berkeley.edu/~jwxiao/}}}. 

Section \ref{subsection_Run time comparison of SRCH and DPSTRF.f} compares the run times of {\tt DPSTRF.f} and {\tt SRCH} to obtain a partial Cholesky factorization. Section \ref{subsection_Approximation effectiveness comparison} compares the approximation effectiveness of low-rank approximations computed by {\tt DPSTRF.f} and {\tt SRCH}. Section \ref{subsection_Test on Kahan matrix} compares {\tt DPSTRF.f} and {\tt SRCH} on a pathological matrix. Section \ref{subsection_Cholesky factorization with side information (CSI)} and section \ref{subsection_Gaussian process} compare {\tt SRCH} with other pivoted Cholesky factorization based algorithms in a prediction problem and Gaussian process, respectively. All experiments are implemented on a laptop with a 2.7 GHz Intel Core i5 CPU and 8GB of RAM.

\subsection{Run time comparison}\label{subsection_Run time comparison of SRCH and DPSTRF.f}
\begin{figure}
\begin{minipage}[t]{0.5\linewidth}
\centering
\includegraphics[width=1.7in]{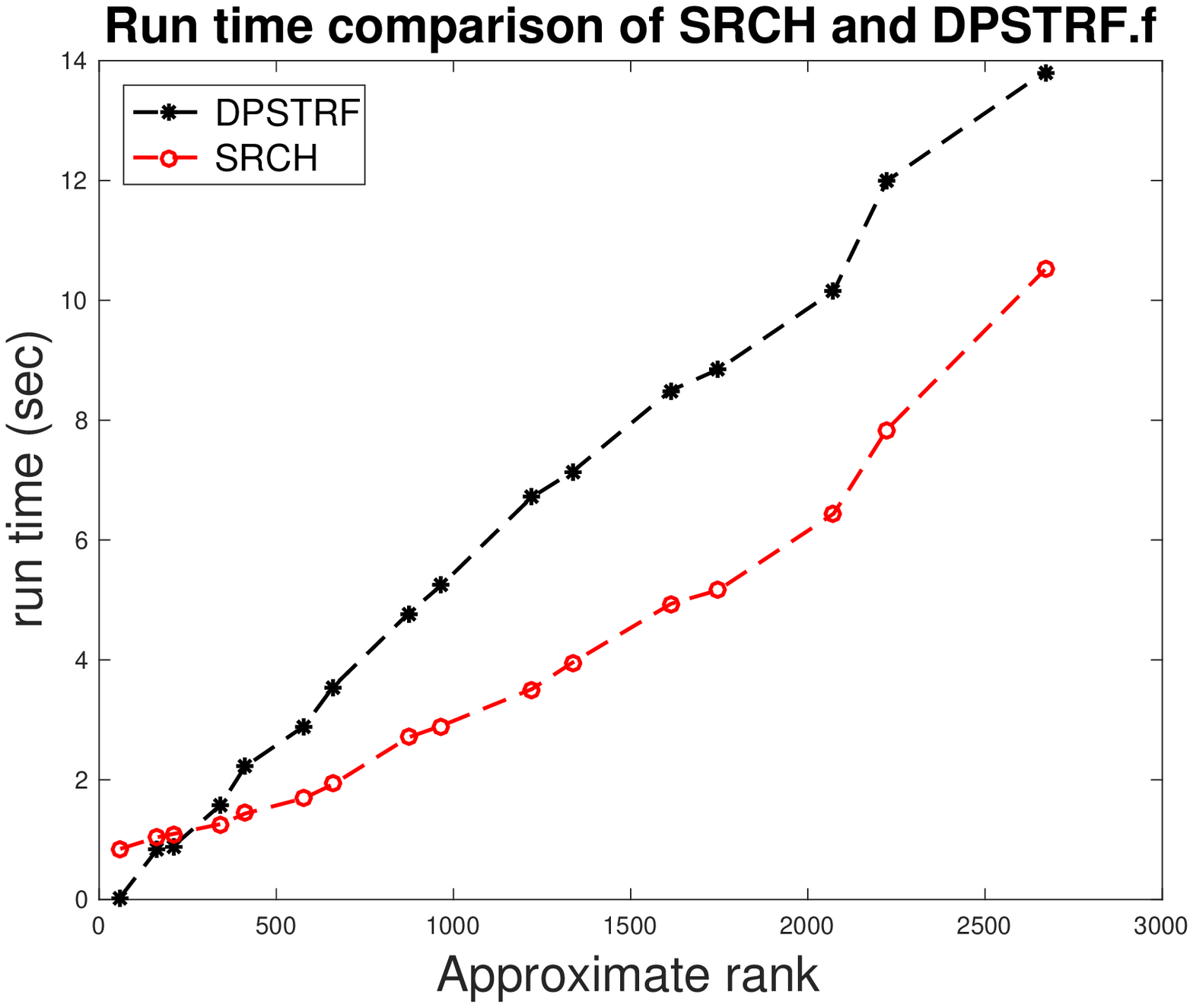}
\caption{Run time comparison.}\label{figure_run_time_srch_dpstrf}
\end{minipage}%
\begin{minipage}[t]{0.5\linewidth}
\centering
\includegraphics[width=1.7in]{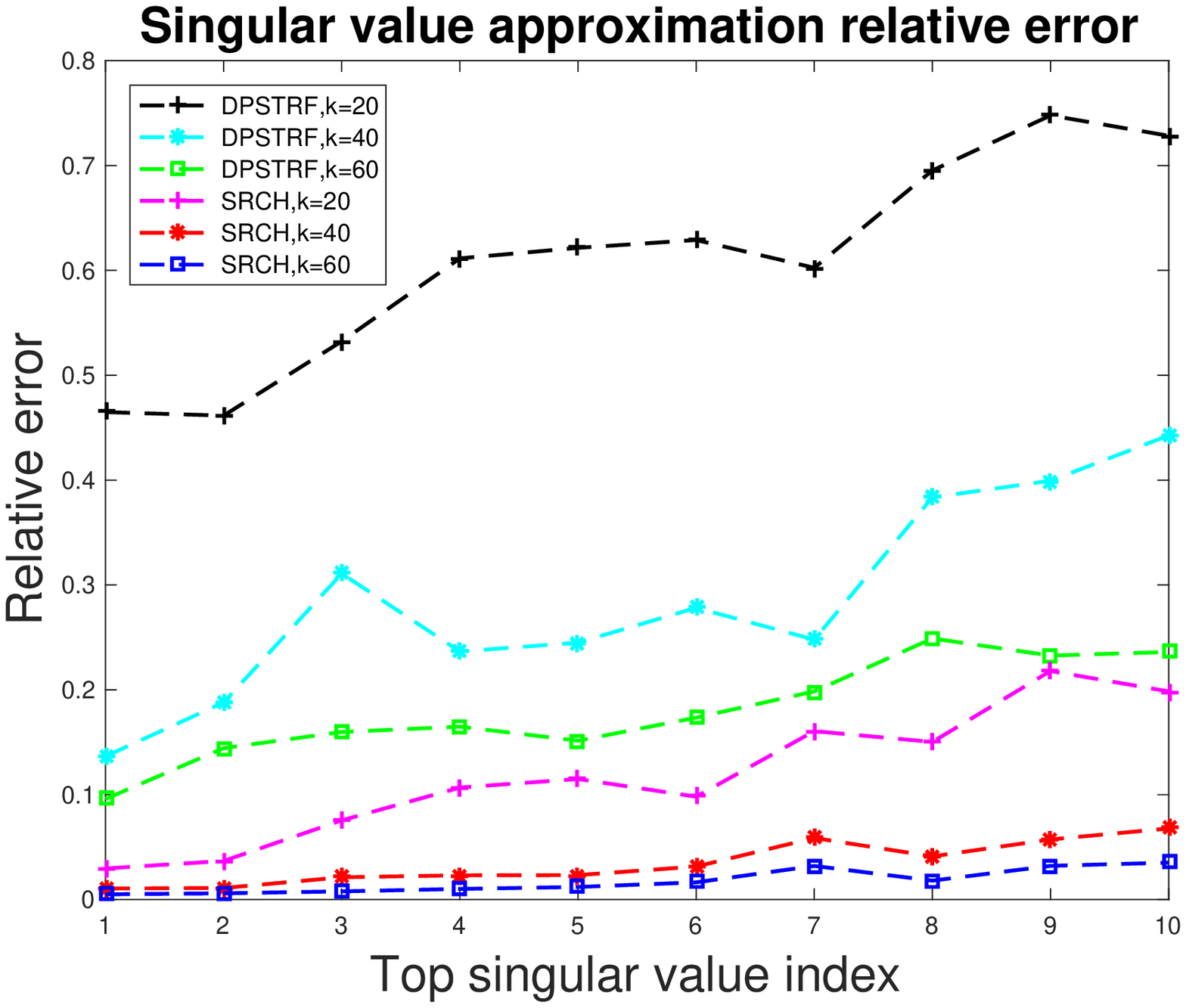}
\caption{Top 10 singular value approximation relative error $\frac{\lambda_{j}(A)-\sigma_{j}^2(L)}{\lambda_{j}(A)}$.}\label{figure_singular_value_approximation_relative_error}
\end{minipage}
\end{figure}

We compare the run times of {\tt DPSTRF.f} and {\tt SRCH} on a kernel matrix of Combined Cycle Power Plant Data (CCPP) with size of $9568\times{}4$. We use RBF kernel $k(x_i,x_j)=\exp{}(-||x_i-x_j||_2^2/2\sigma^2)$ and set $\sigma=1.0$. In {\tt SRCH}, we choose $b=20$ and $p=30$. Run times are in \hyperref[figure_run_time_srch_dpstrf]{Fig.  \ref{figure_run_time_srch_dpstrf}}.

{\tt DPSTRF.f} is faster when approximate rank is small. This is because {\tt SRCH} must compute the random projection $B=\Omega{}A$, which is an overhead. As approximate rank increases, {\tt SRCH} becomes faster than {\tt DPSTRF.f} as predicted. 

\subsection{Approximation effectiveness comparison}\label{subsection_Approximation effectiveness comparison}

For the same kernel matrix used in section \ref{subsection_Run time comparison of SRCH and DPSTRF.f}, we compare the approximation effectiveness of {\tt SRCH} and {\tt DPSTRF.f} with their low-rank approximations.  We choose $k=20, 40, 60$, $b=20$ and $p=30$. \hyperref[figure_singular_value_approximation_relative_error]{Fig. \ref{figure_singular_value_approximation_relative_error}} shows the approximation relative errors of the top $10$ eigenvalues of $A$. Although {\tt SRCH} will be slower than {\tt DPSTRF.f} because of the overhead cost of computing $B=\Omega{}A$, the approximation effectiveness of {\tt SRCH} is much better than {\tt DPSTRF.f}.

If we choose large approximate rank $k$, both {\tt SRCH} and {\tt DPSTRF.f} will provide very high quality low-rank approximations and there will be little difference between {\tt SRCH} and {\tt DPSTRF.f} in relative approximation errors for the leading singular values, but {\tt SRCH} is significantly faster than {\tt DPSTRF.f}, as \hyperref[figure_run_time_srch_dpstrf]{Fig.  \ref{figure_run_time_srch_dpstrf}} suggests.

\subsection{A pathological example for Spectrum-revealing Cholesky factorization}\label{subsection_Test on Kahan matrix}

The Kahan Matrix \cite{kahan1966numerical} is defined as $K_n = S_n C_n$, where 
\footnotesize
$${ S_n ={\textrm{\bf diag}}(1, s, \cdots, s^{n-1}) , C_n = \left(\begin{array}{ccccc}
1 & - c & -c & \cdots & -c \cr
& 1 & -c & \cdots & -c \cr
&  & 1 & \cdots & -c \cr
&  &  & \ddots & \vdots \cr
&  &  &  & 1\end{array}\right)}$$
\normalsize
for $s,c>0$ and $s^2+c^2\leq{}1$. The Kahan matrix is well-known for its peculiar behavior  regarding estimation of conditioning and rank. Indeed, simple diagonal pivoting on the matrix $A\stackrel{def}{=}K_n^T\,K_n$ will fail to produce a quality low-rank approximation. 

We use {\tt DPSTRF.f} and {\tt SRCH} to compute a partial Cholesky factorization of $A$. We choose $c=0.285$, $s = \sqrt{0.9999- c^2}$, $n=130$ and $k=100$. We also set $g=1.5$, $d=b=20$ and $p=25$ in {\tt SRCH}. Let the partial Cholesky factorization be $\Pi^TA\Pi=\left(
\begin{array}{cc}
L_{11} &  \\
L_{21} & I_{n-k} \\
\end{array}
\right)\left(
\begin{array}{cc}
L_{11}^T & L_{21}^T \\
& S \\
\end{array}
\right)$ and we denote $L=\left(
\begin{array}{c}
L_{11}   \\
L_{21}   \\
\end{array}
\right)$. \hyperref[table_Singular value approximation ratio]{Table \ref{table_Singular value approximation ratio}} compares the approximation effectiveness of a few smallest singular values. 

\begin{table}
\begin{center}
\begin{tabular}{|c|c|c|}
\hline
index & DPSTRF & SRCH  \\\hline
96 & 0.8855 &  0.9545 \\\hline  
97 & 0.8739 &  0.9467 \\\hline   
98 & 0.8594 &  0.9370 \\\hline   
99 & 0.8390 &  0.9242 \\\hline   
100 & 0.5820E-08 & 0.9055 \\
\hline
\end{tabular}
\end{center}
\caption{Singular value approximation ratio $\frac{\sigma_j^2(L)}{\lambda_j(A)}$}\label{table_Singular value approximation ratio}
\end{table}

The singular value ratios $\frac{\sigma_j^2(L)}{\lambda_j(A)}$ can never exceed $1$ for any approximation, but we would like them to be close to $1$ for a reliable spectrum-revealing Cholesky factorization. \hyperref[table_Singular value approximation ratio]{Table \ref{table_Singular value approximation ratio}} demonstrates that 
{\tt DPSTRF.f} failed to do so for the index $100$ singular value, whereas {\tt SRCH} has succeeded for all singular values. {\tt SRCH} required $2$ extra swaps to achieve this reliability. The run time of {\tt DPSTRF.f} and {\tt SRCH} are 5.060e-4 seconds and 1.359e-3 seconds respectively. {\tt SRCH} is slower because the Kahan matrix in testing is small in dimension.              

In typical machine learning and other applications of low-rank approximations, the extra swaps of \hyperref[algorithm_swap_strategy]{Algorithm \ref{algorithm_swap_strategy}} are rarely needed. They serve as an insurance policy against occasional mistakes made by \hyperref[algorithm_SRCH]{Algorithm \ref{algorithm_SRCH}}. It is worth noting that the additional running time is negligible.

\subsection{Cholesky factorization with side information (CSI)}\label{subsection_Cholesky factorization with side information (CSI)}

\cite{bach2005predictive} presents an algorithm that exploits side information in the prediction on unlabeled data with low-rank approximations for kernel matrices. To compute a low-rank approximation, this algorithm minimizes the objective function with a greedy strategy to incrementally select representative samples.

In this section we apply {\tt SRCH} on the kernel matrix for a low-rank approximation {\em without} the benefit of side information, and make  predictions on unlabeled data with this low-rank approximation. Details of the prediction formulas can be found in Section $6$ of \cite{bach2005predictive}. 

We compare approximation effectiveness of the low-rank approximation, run time and prediction error on unlabeled data. We compare four methods: CSI decomposition with $40$ look-ahead steps, CSI decomposition with $80$ look-ahead steps, diagonal pivoted Cholesky without look-ahead and {\tt SRCH}. The first three methods are from \cite{bach2005predictive}. We test these four methods on handwritten digit (MNIST). We use RBF kernel $k(x_i,x_j)=\exp{}(-||x_i-x_j||_2^2/2\sigma^2)$ and set $\sigma=1.0$. We set $b=50$ and $p=55$ in {\tt SRCH}.  

\begin{figure}
\begin{minipage}[t]{0.5\linewidth}
\centering
\includegraphics[width=1.7in]{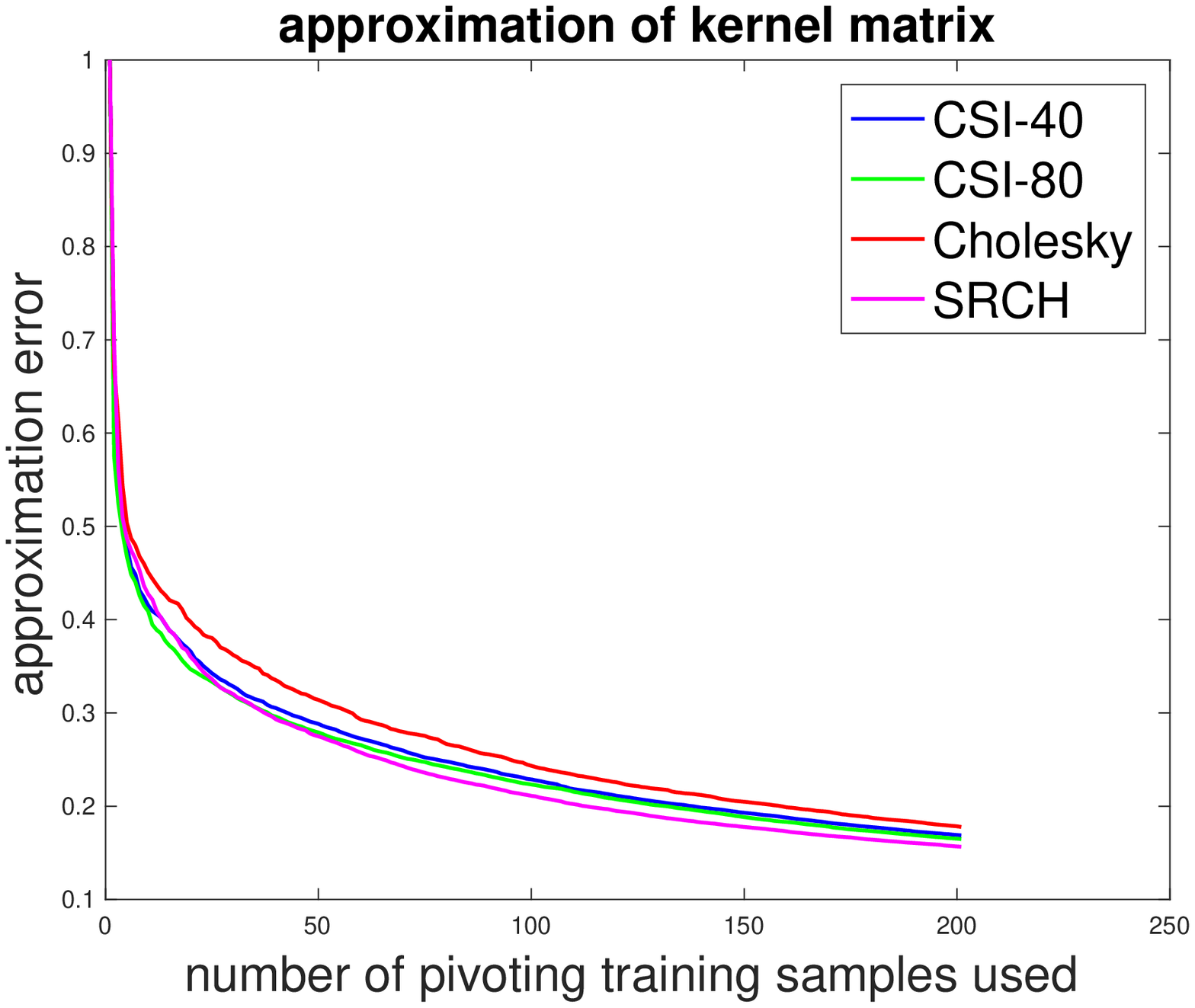}
\caption{Approximation comparison.}\label{figure_mnist approximation of kernel matrix}
\end{minipage}%
\begin{minipage}[t]{0.5\linewidth}
\centering
\includegraphics[width=1.7in]{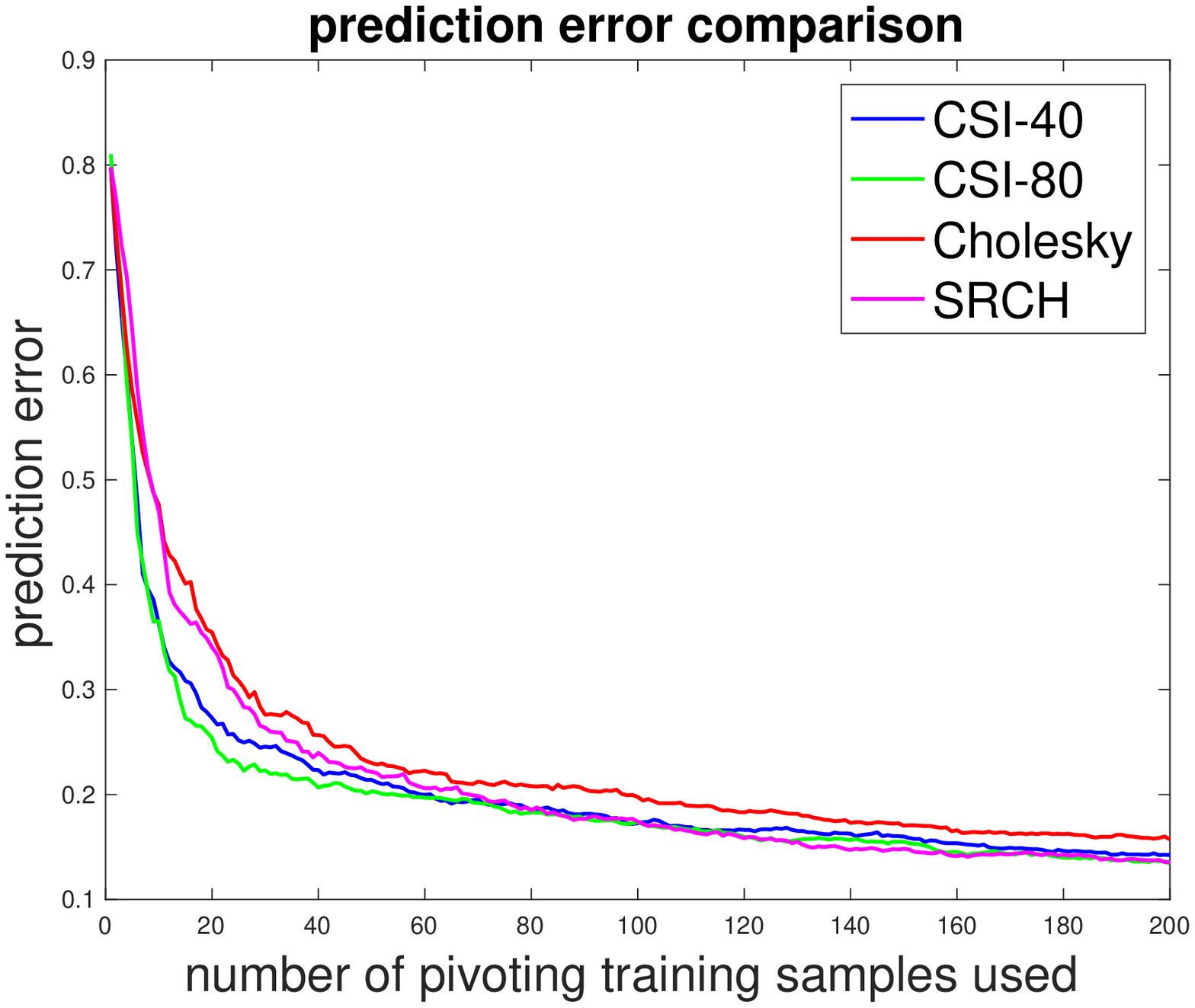}
\caption{Prediction error comparison.}\label{figure_mnist prediction error vs number of pivoting training samples used}
\end{minipage}
\end{figure}

We choose $3000$ training samples, $3000$ testing samples and $k=200$. \hyperref[figure_mnist approximation of kernel matrix]{Fig. \ref{figure_mnist approximation of kernel matrix}} and  \hyperref[figure_mnist prediction error vs number of pivoting training samples used]{Fig. \ref{figure_mnist prediction error vs number of pivoting training samples used}} show approximation effectiveness and prediction accuracy, respectively. There is a slight advantage of {\tt SRCH} on both approximation effect and prediction accuracy. We define the approximation error as $\frac{trace(K-LL^T)}{trace(K)}$.

\begin{figure}
\begin{minipage}[t]{1.0\linewidth}
\centering
\includegraphics[width=1.7in]{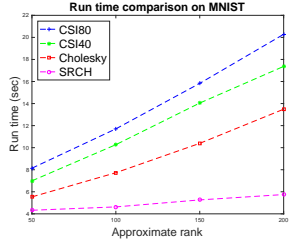}
\caption{Run time comparison.}\label{figure_mnist_running_time_3000}
\end{minipage}
\end{figure}

The more impressive improvement is in run time. \hyperref[figure_mnist_running_time_3000]{Fig. \ref{figure_mnist_running_time_3000}} shows the run time comparison on the kernel matrix $K\in{}\mathbb{R}^{3000\times{}3000}$ for different approximate ranks $k$. {\tt SRCH} is significantly faster than the other three methods. 

\subsection{Gaussian process}\label{subsection_Gaussian process}
Supervised learning is to learn input-output mappings using data. We assume the training data is $X \in \mathbb{R}^{n \times d}$, the target values of $X$ is $y \in \mathbb{R}^{n \times 1}$, the new data is $X^* \in \mathbb{R}^{n^* \times d}$ and the target values of $X^*$ is $y^* \in \mathbb{R}^{n^* \times 1}$. The goal is to predict $y^*$ given $X, y$ and $X^*$ \cite{foster2009stable}. In the Gaussian process, the prediction of $y^{*}$ involves a covariance function $\kappa(x,x^{\prime})$ where $x, x^{\prime} \in \mathbb{R}^d$. The covariance function can be used to construct the $n\times{}n$ covariance matrix $A$ with entries $A_{ij}=\kappa(x_i,x_j)$ where $x_i$ and $x_j$ are rows of $X$ and also the $n^{*}\times{}n$ cross covariance matrix $A^{*}$ where $A^{*}_{ij}=\kappa(x^{*}_i,x_j)$ where $x^{*}_i$ is the $i$th row of $X^{*}$. The prediction $\hat{y}^{*}$ for $y$ is given by the Gaussian process equation \cite{williams2006gaussian}
\begin{equation}\label{equation_Gaussian process equation}
\hat{y}^*=A^*(\lambda I+A)^{-1}y
\end{equation}
where $\lambda{}$ is a regularization parameter. It is not practical to solve equation (\ref{equation_Gaussian process equation}) with large $n$ since the number of flops required is $O(n^3)$. Therefore for large $n$ it is useful to develop approximate solutions to equation (\ref{equation_Gaussian process equation}). We use the approximation formulas in \cite{foster2009stable} in this experiment.

\begin{figure}
\begin{minipage}[t]{0.5\linewidth}
\centering
\includegraphics[width=1.7in]{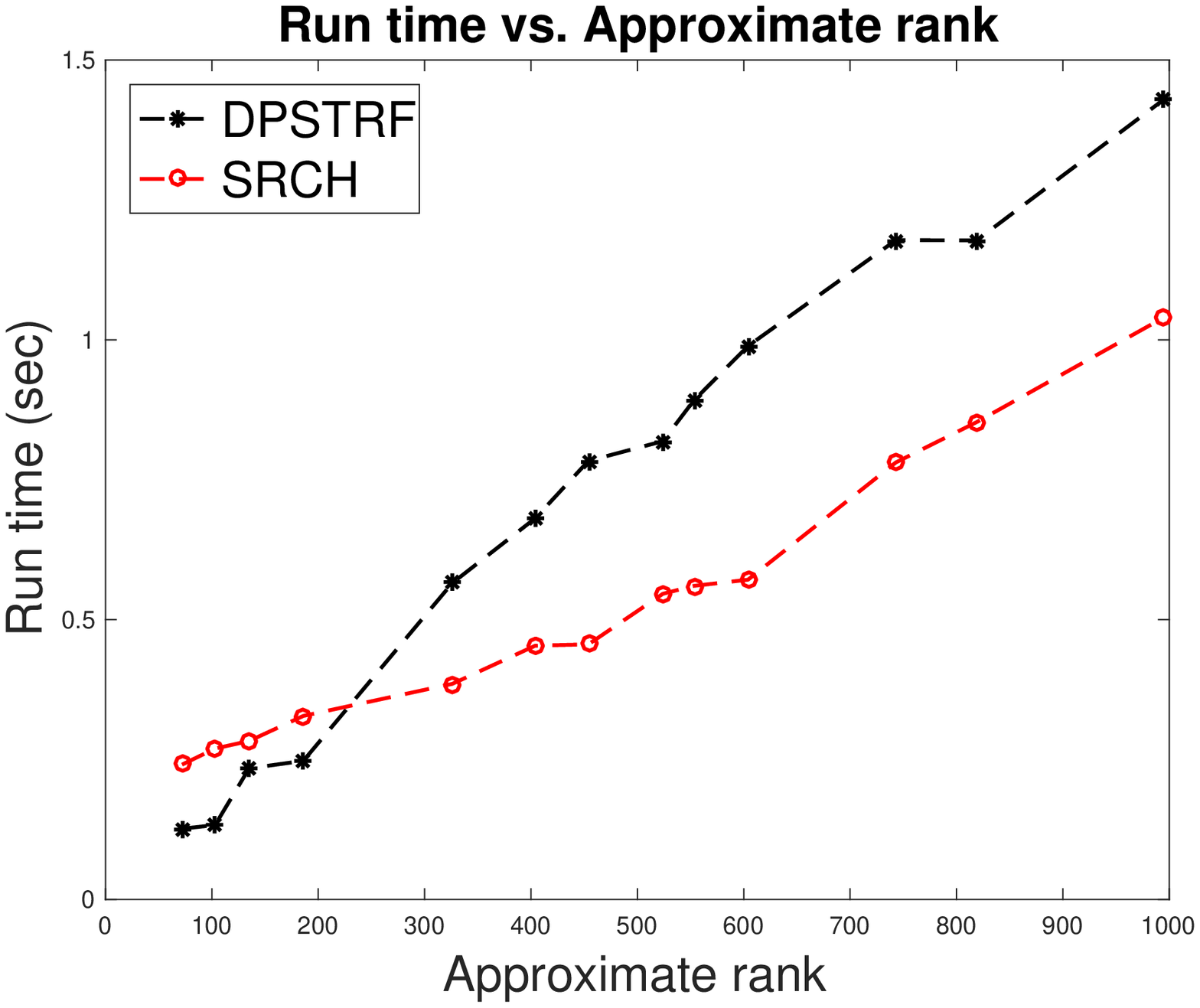}
\caption{Run time comparison.}\label{figure_runningtime_gp_ccpp}
\end{minipage}%
\begin{minipage}[t]{0.5\linewidth}
\centering
\includegraphics[width=1.7in]{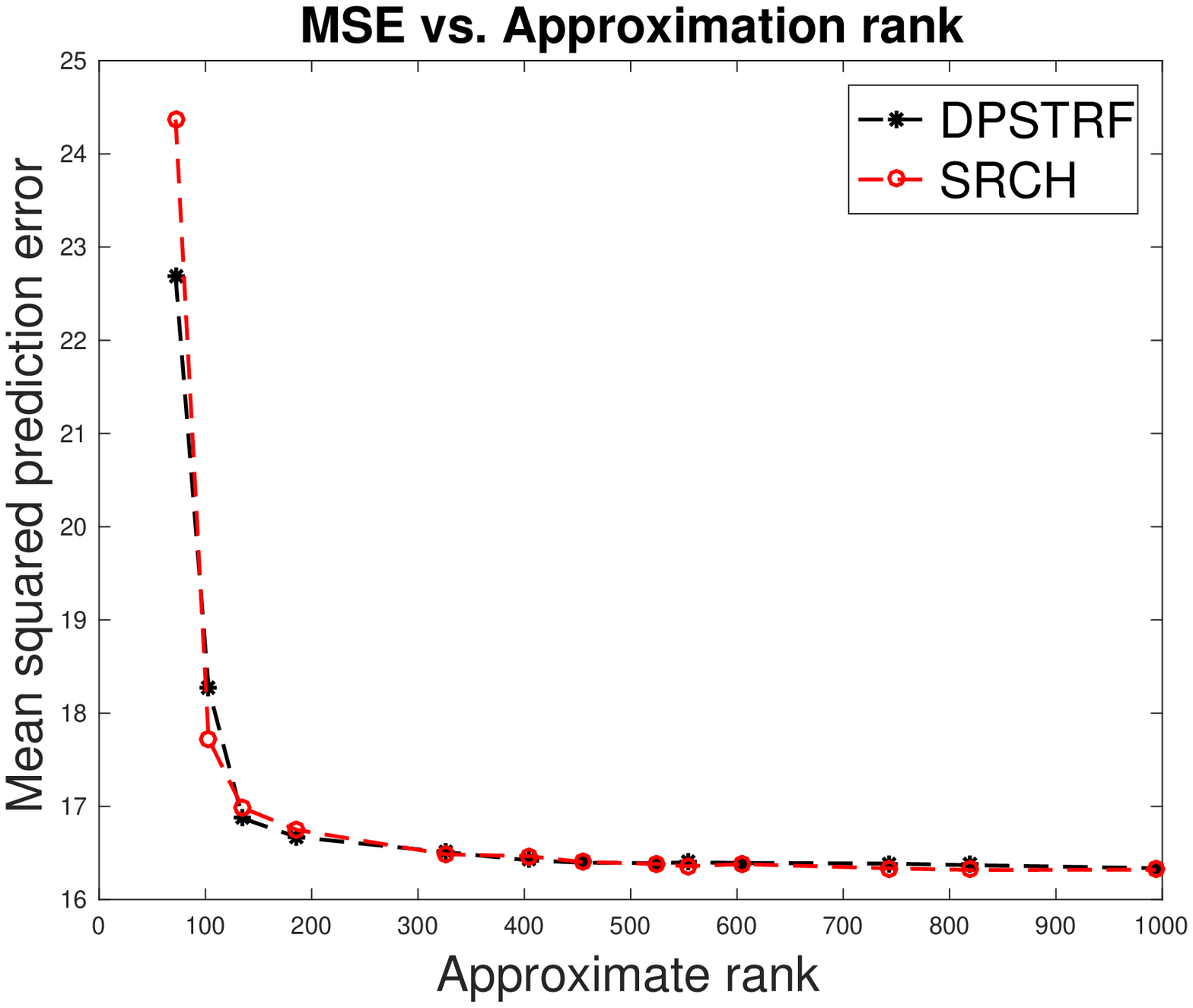}
\caption{Prediction error comparison.}\label{figure_mse_gp_ccpp}
\end{minipage}
\end{figure}

We compute the Gaussian process on CCPP dataset with {\tt DPSTRF.f} and {\tt SRCH}. The training data matrix is of size $5000\times{}4$ and the testing data matrix is of size $4568\times{}4$. The covariance function is $\kappa(x_i,x_j)=\exp{}(-||x_i-x_j||_2^2/2\sigma^2)$. We set $\lambda=5e-5$ and $\sigma=2.0$. In {\tt SRCH}, we set $b=20$ and $p=25$. We compare the run time and mean squared prediction error in \hyperref[figure_runningtime_gp_ccpp]{Fig. \ref{figure_runningtime_gp_ccpp}} and \hyperref[figure_mse_gp_ccpp]{Fig. \ref{figure_mse_gp_ccpp}}, respectively.

In \hyperref[figure_runningtime_gp_ccpp]{Fig. \ref{figure_runningtime_gp_ccpp}}, the run times of {\tt DPSTRF.f} and {\tt SRCH} intersect at around $250$. {\tt SRCH} out-performs {\tt DPSTRF.f} for larger approximate ranks. \hyperref[figure_mse_gp_ccpp]{Fig. \ref{figure_mse_gp_ccpp}} demonstrates that while {\tt DPSTRF.f} makes better predictions than {\tt SRCH} for smaller approximate ranks $k$, they are not the range in which the best predictions are made. For larger $k$, {\tt SRCH} provides slightly smaller prediction error than {\tt DPSTRF.f}, as suggested in \hyperref[figure_runningtime_gp_ccpp]{Fig. \ref{figure_runningtime_gp_ccpp}}.

\section{Conclusion}\label{section_Conclusion}
In this work, we introduced spectrum-revealing Cholesky factorization (SRCH), a variant of the classical Cholesky factorization, for reliable low-rank matrix approximations of kernel matrices. We developed approximation error bounds as well as singular value approximation lower bounds for SRCH. We also developed an efficient and effective randomized algorithm for computing SRCH and demonstrated its efficiency and reliability against other Cholesky factorization based kernel methods on machine learning problems.

% \section*{Acknowledgment}

% \bibliographystyle{plain}
% \bibliography{main}
% \bibliographystyle{unsrt}
\bibliographystyle{IEEEtran}
\bibliography{main.bib} 
\end{document}